\documentclass[11pt]{elsarticle}
\usepackage{amsfonts,amsmath,amsthm}

\usepackage{longtable}

\newcommand{\Z}{\mathbb{Z}}

\newcommand{\abs}[1]{\lvert#1\rvert}
\renewcommand{\phi}{\varphi}

\newcommand{\on}[1]{\operatorname{#1}}
\newcommand{\GF}{\on{GF}}
\newcommand{\GR}{\on{GR}}

\newcommand{\BIGOP}[1]{\mathop{\mathchoice%
{\raise-0.22em\hbox{\huge $#1$}}%
{\raise-0.05em\hbox{\Large $#1$}}{\hbox{\large $#1$}}{#1}}}
\newcommand{\bigtimes}{\BIGOP{\times}}

\newtheorem{theorem}{Theorem}
\newtheorem{corollary}[theorem]{Corollary}

\newtheorem{lemma}[theorem]{Lemma}

\theoremstyle{definition}
\newtheorem{definition}[theorem]{Definition}
\newtheorem{example}[theorem]{Example}

\theoremstyle{remark}
\newtheorem{remark}[theorem]{Remark}

\newcommand{\R}{\mathbb{R}}
\newcommand{\C}{\mathbb{C}}
\newcommand{\Q}{\mathbb{Q}}

\DeclareMathOperator{\soc}{Soc}
\DeclareMathOperator{\rad}{Rad}

\begin{document}

\begin{frontmatter}

\title{Properties of Codes with Two Homogeneous Weights}


%
%

\author[ucd]{Eimear Byrne\fnref{fn1}}
\author[ub]{Michael Kiermaier\fnref{fn2}}
\author[ucd]{Alison Sneyd\fnref{fn1}}

\address[ucd]{%
  School of Mathematical Sciences and the Claude Shannon Institute,\\
  University College Dublin, Ireland}

\address[ub]{%
  Mathematisches Institut,\\
  Universit\"{a}t Bayreuth, Germany}

\fntext[fn1]{Work supported by Science Foundation of Ireland grant 08/RFP/MTH1181}
\fntext[fn2]{Work supported by Deutsche Forschungsgemeinschaft grant WA 1666 4/2}
  
  \begin{abstract}
   Delsarte showed that for any projective linear code over a finite field $\GF(p^r)$ 
   with two nonzero Hamming weights $w_1<w_2$ there exist positive integers $u$ and $s$ such that $ w_1 = p^s u$ and 
   $w_2 = p^s(u+1)$. Moreover, he showed that the additive group of such a code has a strongly regular Cayley graph. 
   Here we show that for any regular projective linear code $C$ over a finite 
   Frobenius ring with two integral nonzero homogeneous 
   weights $w_1<w_2$ there is a positive integer $d$, a divisor of $\abs{C}$, and positive integer $u$ such that
   $w_1 = d u$ and $w_2 = d(u+1)$. This gives a new proof of the known result that any such code yields a 
   strongly regular graph.
   We apply these results to existence questions on two-weight codes.
\end{abstract}   



\end{frontmatter}

\section{Introduction}

        The homogeneous weight 
        has been studied extensively in the context of ring-linear coding. The reader is invited to check 
        the references (cf.~\cite{ch,bruce,hn}) for its different forms. 
        In this paper we will use the definition given in
        \cite{bruce}. 
        The homogeneous weight can in some sense be viewed as a generalization of the Hamming weight; in 
        fact it coincides with the Hamming weight when the underlying ring is a finite field and 
        is the Lee weight when the underlying ring is $\mathbb{Z}_4$. 
        
        Many of the classical results for codes over finite 
        fields for the Hamming weight have 
        corresponding homogeneous weight versions for codes over finite rings. In particular,
        in \cite{bgh} it was shown that strongly regular 
        graphs can be constructed from codes over finite Frobenius rings with exactly two 
        nonzero homogeneous weights. 

        In this paper, we examine properties of the parameters of a regular projective 
        two-weight code over a finite Frobenius ring.
        In Section 3
        we obtain the analogue of a result of Delsarte \cite[Corollary 2]{D72}, namely that 
        such a code $C$ must have nonzero homogenous weights of the form $w_1 = du$ and $w_2 = 
        d(u+1)$ (after scaling by $\abs{R^\times}$), where $d$ is a divisor of $\abs{C}$ and 
        $u$ is a positive integer. This emerges by an analysis of the eigenvalues of the Cayley 
        graph $\Gamma$ of $C$, which turns out to have exactly three distinct eigenvalues and is 
        therefore strongly regular.
        
        We then derive a number of constraints involving the parameters of a strongly regular graph 
        srg$(N,k,\lambda,\mu)$ that can be 
        constructed from a two-weight code and use these to examine the tables of feasible parameter sets of non-trivial
        strongly regular graphs from \cite{b} to see which might arise as Cayley graphs defined from two-weight codes.
        The vast majority of these parameter sets could be eliminated without reference to the underlying ring.  
        For the cases that remained after this, we searched for codes over $R$ of length $n$ determined by the size of the unit group of $R$
        and the graph parameters. An exhaustive search using the techniques described here is feasible only
        for rings of order not divisible by a fourth power (which are then all direct sums of chain rings).
        The open cases that remain are listed in Section 4 of the paper.
                
        In the last section, we show that the known constructions of linear two-weight codes 
        over a chain ring have no linear codes preimages under the Gray isometry.

\section{Preliminaries}
 \subsection{Finite Rings and Homogeneous Weights}
         We recall some properties of finite rings that meet our purposes, many of which are 
         discussed in \cite{h}. 
         An extensive treatment of ring theory can be read in \cite{lam}. See also \cite{lam1,Nechaev-FiniteRings-2008}.
         For a finite ring $R$, we denote by $\hat{R}:=$ Hom$_{\Z}(R,\C^{\times})$,
         the group of additive characters of $R$. $\hat{R}$
         is an $R$-$R$ bi-module according to the relations
         $$^r{\chi}(x) = \chi(rx),\:\:\: \chi^r(x) = \chi(xr)$$ 
         for all $x,r \in R, \chi \in \hat{R}$. 
         A character $\chi$ is called left (resp. right) generating
         if given any $\phi \in \hat{R}$ there is some $r \in R$ satisfying 
         $\phi = {^r}{\chi}$ (resp. $\phi = {\chi}^r$). 
         The next result gives a characterization of finite Frobenius rings.    
         \begin{theorem}
         Let $R$ be a finite ring. The following are equivalent.
         \begin{enumerate}
         \item
         $R$ is a Frobenius ring 
         \item
         $\soc {_R}R$ is a left cyclic module,
         \item
         ${_R}(R/\rad\; R) \simeq  \soc{_R}R$, 
         \item
         $ {_R}R \simeq {_R}\hat{R} $ 
         \end{enumerate}
         \end{theorem}

         Then $~_R\hat{R}= R \chi$ for some (left) generating character $\chi$. 
         It can be shown that any left generating character is also a right 
         generating character (c.f. \cite{wood97a}).

\vspace{1cm} 
 
 For an arbitrary finite ring, the homogeneous weight is defined as follows \cite{bruce}. See also \cite{hn}, for a slightly different definition.  

  \begin{definition}
  Let $R$ be a finite ring. A weight $w:R\longrightarrow \Q$ is
   \emph{(left) homogeneous}, if $w(0)=0$ and
    \begin{enumerate}
  \item If $Rx=Ry$ then $w(x)= w(y)$ for all $x,y\in R$.
  \item There exists a real number $\gamma$ (independent of $R$) 
  such that
    \begin{equation*}
  \sum_{y\in Rx}w(y) \; =\; \gamma \, \abs{Rx}\qquad\text{for all $x\in
    R\setminus \{0\}$}.
\end{equation*}
\end{enumerate}
\end{definition}

Right homogeneous weights are defined similarly.


 \begin{example}
  On every finite field $\GF(q)$ the Hamming weight
    is a homogeneous weight of average value $\gamma=\frac{q-1}{q}$.
 \end{example}

\begin{example}
On the ring $\Z_{pq}, p,q$ prime, a homogeneous weight with average value $\gamma = 1$ is given by
 \begin{equation*}
      w: R\longrightarrow {\mathbb R},\quad x \mapsto
      \left\{\begin{array}{cl}
                     0& \mbox{if } x=0,\\
        \displaystyle{\frac{p}{p-1}} & \mbox{if } x \in p\Z_{pq} ,\\
        \displaystyle{\frac{q}{q-1}}& \mbox{if } x \in q\Z_{pq},\\
        \displaystyle{\frac{pq-p-q}{(p-1)(q-1)}} & {\mbox{ otherwise.}}       
      \end{array}\right.
    \end{equation*}
\end{example}

In fact the homogeneous weight is unique up to choice of $\gamma$ on any finite ring. The definition used in \cite{ch,bruce} uses the notion of a {\rm M\"{o}bius function} on the set of left principal ideals of $R$, given the partial order induced by set inclusion. It an integer-valued function implicitly defined by
\begin{eqnarray*}
\mu(Rx,Rx) & = & 1 \;\;\mbox{for all $x\in R$},\\
\mu(Ry,Rx) & = & 0 \;\;\mbox{if $Ry \not \leq Rx$, and}\\
\sum_{Ry \leq Rz \leq Rx}\mu(Rz,Rx) & = &  0,\;\;\mbox{if $Ry < Rx$.}
\end{eqnarray*}

\begin{definition}\label{defmobhom}
   Let $\mu$ denote the M\"{o}bius function on the partially ordered set of left principal ideals of $R$.
   For any $r \in R$ the homogeneous weight of $r$ is given by
   $$w(r) \;=\; \gamma\left(1-\frac{\mu(0,Rx)}{\abs{R^\times x}}\right),$$
   where $\gamma$ is a real constant. We say that $w$ is the {\em normalized} 
   homogeneous weight for the case $\gamma = 1$.
\end{definition}

Clearly, $w(r) \in \Q$ whenever $\gamma \in \Q$. In particular, if $\gamma = \abs{R^\times}$ then $w$ is 
an integer-valued function, since $\abs{R^{\times}} = \abs{R^\times x}\abs{{\rm Stab}_{R^{\times}}(x)}$ for all $x \in R$. 

\begin{example}
On a local Frobenius ring $R$ with residue field $\GF(q)$, 
we have $\mu(0,Rx)=-1$ for $Rx = \soc R$ and $\mu(0,Rx)=0$ for $x \in R \backslash \soc R$.
    The homogeneous weight is given by
    \begin{equation*}
      w: R\longrightarrow {\mathbb R},\quad x \mapsto
      \left\{\begin{array}{cl}
        0       & \mbox{if } x=0,\\
        q       & \mbox{if } x\in \soc(R),\; x\neq 0,\\
        q-1     & \mbox{if }  \mbox{otherwise},
      \end{array}\right.
    \end{equation*}
where we choose $\gamma=q-1$.
\end{example}

\begin{lemma}\label{lemsum}
    Let $R$ and $S$ be finite rings. Let $w_R$ and $w_S$ be homogeneous weights on $R$ and $S$ respectively. Then the homogeneous weight of $(r,s) \in R \oplus S$, denoted $w_{R\oplus S}(r,s)$, satisfies
    $$1-\frac{w_{R\oplus S}(r,s)}{\gamma_{R\oplus S}} =  \left(1-\frac{w_R(r)}{\gamma_R} \right) \left(1-\frac{w_S(s)}{\gamma_S} \right),$$
    for real numbers $\gamma_{R \oplus S}, \gamma_R$ and $\gamma_S$.
\end{lemma}

{\bf Proof: }
The M\"{o}bius function on the partially ordered set of left principal ideals of $R\oplus S$
is given by 
$$\mu_{R \oplus S}((Rr,Ss),(Rr',Ss'))=\mu_R(Rr,Rr')\mu_S(Rs,Rs'),$$
where $\mu_R$ and $\mu_S$ are the M\"{o}bius functions on the partially ordered sets of left principal ideals of $R$ and $S$, respectively. The rest follows directly from Definition \ref{defmobhom}.
\qed 

A description of the homogeneous weight in terms of sums of generating characters is given by the following \cite{h}.
\begin{theorem} 
Let $R$ be a finite Frobenius
  ring with generating character $\chi$. 
  Then the homogeneous weights on $R$ are precisely the functions
  \begin{equation*}
w: R\longrightarrow {\mathbb R}, \quad x \mapsto
  \gamma\left(1-\frac{1}{\abs{R^{\times}}}\sum_{u\in R^{\times}}
  \chi(xu)\right)
\end{equation*}
where $\gamma$ is a real number. The normalized homogeneous weight occurs for $\gamma=1$.
\end{theorem}


\subsection{Codes Over Rings}

For the remainder, unless otherwise stated, 
we let $R$ denote a finite Frobenius ring endowed with a homogeneous weight $w$.
We extend $w$ to a weight function on $R^n$ in the obvious way:
$$w: R^n \longrightarrow \R:w(c_1,\ldots,c_n) \mapsto \sum_{i=1}^n w(c_i).$$
We also let $C \leq {_RR^n}$ denote a left linear code. As usual, $I$ and $J$ will denote the real identity matrix and the real all-ones matrix, respectively.

\begin{definition}\label{defproj}
  Let $C$ have $\ell \times n$ generator matrix
  $Y=(y_1|y_2|\dots|y_n)$ over $R$ (we do not assume that the rows of $Y$ are linearly independent over $R$). 
  $C$ is called
  \begin{enumerate}
  \item
     proper if $w(c)=0$ implies $c=0$,\footnote{Note that the homogeneous weight is not positive definite on all finite Frobenius rings; for example in the ring $\GF(2) \times \GF(2)$ we have $w(00)=w(11)=0$.} for all $c \in C$.
  \item regular if $\{x \cdot y_i : x \in R^\ell\} = R$ for each 
  $i \in \{1,\dots,n\}$,
  \item projective if $y_iR\neq y_jR$ for any pair of distinct
    coordinates $i,j \in \{1,\dots,n\}$,
  \end{enumerate}
\end{definition}

  We remark that these notions are independent of the particular choice of $Y$. 

\subsection{Strongly Regular Graphs}


We recall some elementary facts about strongly regular graphs. There are many texts on the subject. The reader is referred to \cite{god}, for example, for further details.

\begin{definition}
A graph $G$ on $N$ vertices is called strongly regular with parameters $(N,k,\lambda,\mu)$ if 
\begin{enumerate}
\item
   $G$ is regular of degree $k$,
\item   
   every pair of adjacent vertices has exactly $\lambda$ common neighbours,
\item
   every pair of non-adjacent vertices has exactly $\mu$ common neighbours.    
\end{enumerate} 
\end{definition}

A strongly regular graph $G$ with parameters $(N,k,\lambda,\mu)$ is called {\em trivial} or {\em imprimitive} if either $G$ or its complement is disconnected. If $G$ is disconnected then $\mu = 0$ and $k = \lambda + 1$, in which case it is the union of some number of complete graphs. If the complement of $G$ is disconnected then $k = \mu$.

Let $A$ be the adjacency matrix of a graph $G$ that is neither null nor complete. Then
$G$ is strongly regular if and only if 
\begin{equation}\label{eqsrg}
AJ=JA=kJ \mbox{ and } A^2-(\lambda-\mu)A - (k-\mu)I = \mu J. 
\end{equation}
An eigenvalue $\rho \neq k$ is called {\em restricted} if it has eigenvector orthogonal to ${\bf 1}$.  
It is well known that a connected regular graph of degree $k$ is strongly regular if and only if it has exactly two distinct restricted eigenvalues.
If $G$ is strongly regular then 
(\ref{eqsrg}) yields that any restricted eigenvalue $\rho$ of $A$ 
satisfies
\begin{equation}\label{eqeigen1}
    \rho^2 - (\lambda-\mu) \rho - (k-\mu)=0,    
\end{equation}
which immediately gives the following.
\begin{lemma}\label{lemsrg}
   Let $G$ be a strongly regular graph with parameters $(N,k,\lambda,\mu)$.
   Let $G$ have adjacency matrix with eigenvalues $\rho_1,\rho_2,k$. Then
   \begin{equation}
       \mu = k + \rho_1\rho_2 \mbox{ and } \lambda = k+ \rho_1+\rho_2 + \rho_1\rho_2
   \end{equation}
\end{lemma}
In particular, the parameters $\lambda, \mu$ are completely determined by $k,\rho_1,\rho_2$.
It can be deduced from the fact that $A$ has zero trace that if $G$ is non-trivial, unless the $\rho_i$ occur with the same multiplicity, they are integers of opposite sign, say, $k>\rho_2>0>\rho_1$ and further, that $\rho_1<-1$. 




In \cite{D72} it was shown
that every projective code over a finite field with exactly two non-zero Hamming weights (also called a two-weight code) has a strongly regular Cayley graph. 

Given a two-weight code $C$, with non-zero weights $w_1<w_2$, we denote by $\Gamma(C)$ the graph whose vertices are the codewords of $C$ and whose edges are pairs of vertices $(c,c')$ such that $w(c-c')=w_1$. $\Gamma(C)$ is the Cayley graph of the set of codewords of weight $w_1$ in $C$. 
If $C$ is a two-weight code, we say that $C$ is imprimitive if $\Gamma(C)$ is trivial. Otherwise we say that $C$ is primitive.

\section{Main Results}

We now determine relations between the eigenvalues of $\Gamma(C)$ and the weights of a two-weight code $C$. 

\begin{definition}
The distance matrix of $C$ is the $\abs{C} \times \abs{C}$ matrix $D$ with rows and columns indexed by
the elements of $C$ and whose $(u,v)$-th entry is ${D}_{uv} = w(u-v)$.
\end{definition}

The following is an extension of \cite[Theorem 1]{D72}. The proof makes use of the character description of the homogeneous weight in the distance matrix of the code.

\begin{theorem}\label{thDJ}
    Let $C < {_R}R^n$ be regular and projective. Let $J$ denote the $\abs{C} \times \abs{C}$ all-ones matrix. Then 
    \begin{enumerate}
    \item[\rm{(i)}]\label{theq1}
    $DJ = \gamma n \abs{C} J$ and
    \item[\rm{(ii)}]\label{theq2}
    $\displaystyle{D^2 + \frac{\abs{C}\gamma}{\abs{R^\times}} D=n \gamma^2 \abs{C}\left( \frac{1}{\abs{R^\times}} +n \right)J}$.
    \end{enumerate}
\end{theorem}

{\bf Proof:}
Let $X$ be the $\abs{C} \times n$ matrix whose rows are the codewords of $C$ and let $B$ be the 
$\abs{C} \times \abs{R^{\times}}n$ matrix $B = [ X\lambda]_{\lambda \in R^\times}.$ 
Each coefficient of $B$ is indexed by the symbols $c,(j,\lambda)$ for each $c \in C,\lambda \in R^\times, j \in \{1,\ldots,n\}$, and we write $B_{c,(j,\lambda)}= c_j\lambda$.
We write $\chi_B$ to denote the complex $\abs{C} \times \abs{R^{\times}}n$ matrix whose components
satisfy $(\chi_B)_{c,(\lambda,i)}:=\chi( c_i\lambda)$, i.e. $\chi_B$ is obtained from $B$ by applying the character $\chi$ to the coefficients of $B$. We write $\chi_B^*$ to denote the adjoint of $\chi_B$.
Let $\pi_i: R^n \longrightarrow R$ denote the projection onto the $i$th coordinate. Since $C$ is regular, $\pi_i(C) = R$ and so
\begin{eqnarray*}
(\chi_B^* {\bf 1} )_{\lambda,i} &=& \sum_{c\in C} \chi(-c_i\lambda) =  \sum_{c\in C} \chi(-\pi_i(c)\lambda)
 = \sum_{r \in R} \chi(-r\lambda)\abs{\{c \in C : \pi_i(c)=r\}}\\
 & = & \abs{\ker \pi_i \cap C} \sum_{r \in R} \chi(r)=0 
\end{eqnarray*} 
for each $\lambda\in R^{\times},i \in \{1,\ldots,n\},$
and hence 
\begin{equation}\label{eq1}
\chi_B^* J = 0.
\end{equation}

Let $\pi_{ij}: R^{n} \longrightarrow R^2$ be the projection of a word in $R^{n}$ onto the pair of coordinates indexed by $i$ and $j$. 
For each $\lambda,\mu \in R^{\times}$, define 
$\theta_{\lambda,\mu}:R^2 \longrightarrow R: (a,b) \mapsto b \mu -a\lambda$. 
Let $\Lambda=\Lambda_{(i,\lambda),(j,\mu)}:=\theta_{\lambda,\mu} \circ \pi_{ij}$. Then
\begin{eqnarray*}
(\chi_B^* \chi_B)_{(i,\lambda),(j,\mu)} 
& = & \sum_{c \in C}\chi\left( c_j\mu- c_i\lambda \right)
                                        = \sum_{r \in \Lambda(C)} \chi(r)\abs{\{ c \in C : \Lambda(c) = r \}}\\
  & = &  \abs{\ker \Lambda \cap C} \sum_{r \in \Lambda(C)} \chi(r)                                                                        
  = \left\{ \begin{array}{cl} 
                 \abs{\ker \Lambda \cap C} & \mbox{if }\Lambda(C) = \{0\},\\
                                    0 & \mbox{otherwise.}
                                                  \end{array} \right.
\end{eqnarray*} 
Clearly, if $\Lambda(C)=\{0\}$ then $\abs{\ker \Lambda \cap C}=\abs{C}$.
Let $Y$ be an $\ell \times n$ generator matrix for $C$.
Now $\Lambda(c)=0$ for all $c \in C$ if and only if there exists a pair of coordinate positions $i,j$ 
such that $c_i\lambda = c_j \mu$ for all $c \in C$, which holds if and only if 
$x\cdot (y_i\lambda - y_j \mu) = 0$ for all $x \in R^\ell$, 
in which case $y_iR = y_jR$. By the assumption that $C$ is projective, we deduce that $\Lambda$ is identically zero on $C$ only if $i=j$. For $i=j$, $\Lambda(c)=0$ for all $c \in C$ if and only if $\pi_i(c)(\lambda - \mu) =0$ for all $c \in C$. Since $C$ is regular, $\pi_i(C)=R$, so this holds if and only if $\lambda=\mu$.   
It follows that 
\begin{equation}\label{eq2}
\chi_B^* \chi_B = \abs{C} I.
\end{equation}

We can also relate $\chi_B$ to the distance matrix of $C$:
\begin{eqnarray*}
(\chi_B \chi_B^*)_{a,b} = \sum_{i=1}^n \sum_{\lambda \in R^\times} \chi( a_i \lambda-b_i \lambda)
                                       & = & \sum_{i=1}^n \abs{R^\times} \left(1 - \frac{1}{\gamma}w(a_i-b_i)\right)\\
                                       & = & \abs{R^\times} \left(n -\frac{1}{\gamma} w(a-b)\right),
\end{eqnarray*}  
which gives 
\begin{equation}\label{eq3}
\chi_B \chi_B^* = \abs{R^\times} \left(n J -\frac{1}{\gamma} D\right)
\end{equation}

Combining (\ref{eq1}) and (\ref{eq3}) we obtain $DJ = \gamma n\abs{C}J$.

Using (\ref{eq1}), (\ref{eq2}) and (\ref{eq3}) we get
\begin{eqnarray*}
    (\chi_B^* \chi_B) \chi_B^* = \abs{C} \chi_B^* =  \chi_B^* (\chi_B \chi_B^*)= \abs{R^\times}\chi_B^*\left(J-\frac{1}{\gamma}D\right) = -\frac{\abs{R^\times}}{\gamma} \chi_B^*D
\end{eqnarray*}
and so $\displaystyle{ \abs{C} \chi_B^* +\frac{\abs{R^\times}}{\gamma} \chi_B^*D=0.}$ Therefore,
 \begin{eqnarray*}
     0 &=& \chi_B  \chi_B^* \left(\abs{C}I+\frac{\abs{R^\times}}{\gamma}D\right)
       = \left(nJ - \frac{1}{\gamma} D\right)\left(\abs{C}I+\frac{\abs{R^\times}}{\gamma}D\right) \\
  \implies 0 &=& n \gamma^2 \abs{C}\left( \frac{1}{\abs{R^\times}} +n \right)J - \frac{\abs{C}\gamma}
  {\abs{R^\times}} D - D^2
\end{eqnarray*}
\qed

The first part of the following was proved in \cite[Th. 5.5]{bgh}, using different techniques.

\begin{corollary}\label{coreigen}
   Let $C$ be a proper, regular, projective two-weight code with nonzero weights $w_1<w_2$. Then 
   $\Gamma:=\Gamma(C)$ is strongly regular and the eigenvalues $k,\rho_1,\rho_2$ of the 
   adjacency matrix of $\Gamma$ satisfy
   \begin{enumerate}
   \item[\rm{(i)}]\label{coreq1}
      $ (w_2-w_1)k =  w_2(\abs{C}-1) - \gamma n \abs{C}$
   \item[\rm{(ii)}]\label{coreq2}
      $(w_2-w_1)\rho_1 = - w_2 $ 
   \item[\rm{(iii)}]\label{coreq3}
      $\displaystyle{(w_2-w_1)\rho_2 =  - w_2 + \frac{ \gamma \abs{C}}{\abs{R^\times}}}$      
   \end{enumerate}
\end{corollary}

{\bf Proof:}
   Since $C$ is proper, the adjacency matrix $A$ of $\Gamma$ satisfies 
   \begin{equation}\label{eqad}
       (w_2-w_1)A = w_2(J-I) - D.
   \end{equation}
   $A,D,J$ are real symmetric commuting matrices and can thus be simultaneously diagonalized by 
   an orthogonal matrix. Applying Theorem \ref{thDJ}, (i), 
   we observe that ${\bf 1}$ is an eigenvector of $A$ with eigenvalue $k$ satisfying 
   (i), above.
   Any eigenvector $e$ of $A$ orthogonal to ${\bf 1}$ satisfies
   $$(w_2-w_1) A e = (w_2-w_1) \rho e = w_2(J-I)e -De= -(w_2+\theta)e $$
   where $\rho,\theta$ are the associated eigenvalues for $A$ and $D$, respectively. 
   From Theorem \ref{thDJ}, (ii)
   , we have $(D-n\gamma \abs{C} I)(D)(D+ \frac{\gamma \abs{C}}{\abs{R^\times}} I)=0$, and hence $D$ has exactly two 
   eigenvalues $\theta_1 = 0$ and $\theta_2 = - \frac{\gamma \abs{C}}{\abs{R^\times}}$ corresponding to eigenvectors 
   orthogonal to ${\bf 1}$. 
   That $A^2$ is the required linear combination of $A,J,I$ follows from (\ref{eqad}) and Theorem \ref{thDJ}, (ii).    
\qed

\begin{corollary}\label{corweight}
   Let $C$ be a proper, regular, projective two-weight code with nonzero weights $w_1<w_2$. 
   Let the adjacency matrix of $\Gamma(C)$ have restricted eigenvalues $\rho_1<\rho_2$. 
   Then 
   \begin{enumerate}
   \item[\rm{(i)}] 
      $\rho_2 - \rho_1$ is an integral divisor of $\abs{C}$;
   \item[\rm{(ii)}] 
   $\displaystyle{w_1 = \frac{\gamma \abs{C}(\rho_1 + 1)}{(\rho_1-\rho_2)\abs{R^\times}}\quad\text{and}\quad
   w_2 = \frac{\gamma \abs{C}\rho_1}{(\rho_1-\rho_2)\abs{R^\times}}\text{.}}
   $
   \end{enumerate}
\end{corollary}

{\bf Proof:}
  From Corollary \ref{coreigen}, $(w_2-w_1)(\rho_2-\rho_1) =\frac{\gamma \abs{C}}{\abs{R^\times}}$.
  For $\gamma = \abs{R^\times}$, we get that the integer $\rho_2-\rho_1$ divides $\abs{C}$; 
  since the values $\rho_1$ and $\rho_2$ are independent of this $\gamma$, this must hold true in general.
  Let $d = w_2-w_1$. Now solve for $w_1$ and $w_2$ using the equations $d\rho_1 = - w_2$ 
  and $d = w_2-w_1$.  
\qed

\begin{corollary}\label{corgamma}
   Let $C$ be a proper, regular, projective two-weight code with nonzero weights $w_1<w_2$ where the 
   weight function is computed for $\gamma = \abs{R^\times}$.
   Then there exists a positive integer 
   $d$, a divisor of $\abs{C}$, and positive integer $t$ such 
   that $w_1=d t$ and $w_2 = d(t+1)$. 
\end{corollary}   

{\bf Proof: }
   For $\gamma = \abs{R^\times}$, $w_1,w_2$ are integers, and their difference $d=w_2-w_1$
   is a positive integer dividing $\abs{C}$. Then from Corollary \ref{coreigen}, (ii), 
   and with the same notation, $\rho_1$ is a negative integer less than $-1$. 
   Moreover, from Corollary \ref{corweight}, we get 
   $w_1 = d(-\rho_1-1)$ and $w_2 = d (-\rho_1)$. The result follows with $t=-(\rho_1+1)$. 
\qed

\begin{corollary}\label{cormulti}
   Let $C$ be a proper, regular, projective two-weight code. 
   Let the adjacency matrix of $\Gamma(C)$ have simple eigenvalue $k$ and 
   restricted eigenvalues $\rho_1<\rho_2$.
   Then the multiplicities $m_1$ and $m_2$ of $\rho_1$ and $\rho_2$, respectively, are given by
   $$m_1 = \abs{C} - 1 - n\abs{R^\times} \mbox{ and } m_2 =  n\abs{R^\times}. $$ 
\end{corollary}

{\bf Proof: }
   This follows immediately from Corollary \ref{coreigen} and the equations (cf. \cite{god})
   \begin{equation*}
        m_1 = \frac{(\abs{C}-1)\rho_2 + k}{\rho_2-\rho_1}
       \mbox{ and }
        m_2 = \frac{(\abs{C}-1)\rho_1 + k}{\rho_1-\rho_2}
   \end{equation*}
   where $\Gamma(C)$ is regular of degree $k$.
\qed

We will use the following simple observation.
\begin{lemma}
\label{lma:nrunits}
   Let $R$ have order $\abs{R}=p^r$, $p$ prime, with a minimal right ideal of size $p^s$.
   Then $p^s-1$ divides $\abs{R^\times}$.
\end{lemma}

{\bf Proof: }
Let $I$ be a minimal right ideal of size $p^s$.
From the minimality we get $I = aR$ with $a\in R$.
For $\gamma = \abs{R^\times}$ we have 
\[
w(a) = \frac{\gamma p^s}{p^s-1} = \frac{\abs{R^\times} p^s}{p^s-1} \in \Z\text{.}
\]
So $p^s-1$ divides $\abs{R^\times}$.
\qed

   Now let $C$ be a proper, regular, projective two-weight code of order $\prod_{i=1}^d p_i^{t_i}$ 
   with nonzero weights $w_1<w_2$ where the 
   weight function is computed for $\gamma = \abs{R^\times}$. 
   Let the largest restricted eigenvalue of the adjacency matrix of $\Gamma(C)$
   have multiplicity $m_2$.
    Since $C$ is an $R$-module, each $p_i$ divides $\abs{R}$ and since $C$ is regular, $\abs{R}$ divides 
    $\abs{C}$. Therefore, $R$ has a minimal right ideal $a_iR$ of order 
    $p_i^{s_i}$ for some $s_i \leq t_i$ and for $\gamma=\abs{R^\times}$ we have
    \begin{equation}\label{eqp1}
       \displaystyle{w(a_i) =  \abs{R^\times}\frac{p_i^{s_i}}{p_i^{s_i}-1}}\text{.}
    \end{equation}
    Since $C$ is regular, there is some $c\in C$ with a unit entry in some coordinate.
    Since the nonzero entries of $a_ic$ are unit multiples of $a_i$ in $R$ we deduce that 
    $w(a_ic)=\ell_i w(a_i)$ for some integer $\ell_i \in \{1,\ldots,n\}$. 
    Now $w(a_ic) \in \{w_1,w_2\}$, so there is a $j\in\{1,2\}$ such that $\ell_i$ divides $w_j$, and by (\ref{eqp1}) 
    \begin{equation}\label{eqalcrit}
    \displaystyle{p_i^{s_i} = \frac{w_jn}{w_jn - m_2\ell_i}}\text{.}
    \end{equation}
    Moreover, from Lemma \ref{lma:nrunits}, $\prod_{i=1}^d (p_i^{s_i}-1)$ divides $\frac{m_2}{n}=\abs{R^\times}$.
    
    Therefore, given a parameter set $(N,k,\lambda,\mu)$, we can often eliminate the possibility that the 
    corresponding graph might arise from a two-weight code of length $n$ without consideration of the underlying ring $R$.
    More generally, we have the following result.

\begin{corollary}\label{coralext}
Let $R = \oplus_{i = 1}^d R_i$ satisfy $\abs{R_i} = p_i^{r_i}$ for distinct primes $p_i$.
Let further $p_i^{s_i}$ with $p_i$ prime be the size of a minimal right ideal in $R_i$.
\begin{enumerate}[(a)]
\item $\abs{R^\times}$ is divisible by $\prod_{i = 1}^d (p_i^{s_i} - 1)$.
In particular for any $T\subseteq \{1,\ldots,d\}$ the number
\[
W_T = \abs{R^\times}\left(1-\prod_{i\in T}\left(-\frac{1}{p_i^{s_i}-1}\right)\right)
\]
is integral.
\item Assume that there exists a proper, regular, projective $R$-linear two-weight code of length~$n$ with nonzero weights $w_1 < w_2$ where the 
weight function is computed for $\gamma = \abs{R^\times}$. 
Then there exist $2^d$ non-negative integers $x_T$ with $T\subseteq\{1,\ldots,d\}$ satisfying the $2^d$ conditions
\begin{enumerate}[(i)]
\item $\sum_{T\subseteq{\{1,\ldots,d\}}}x_T = n$ and
\item For all $T\subseteq\{1,\ldots,d\}$, $T\neq\emptyset$:
\[
\sum_{U\subseteq T} \left(\sum_{{V\subseteq\{1,\ldots,d\}}\atop{V\cap T = U}} x_V\right) W_U\in \{w_1,w_2\}\text{.}
\]
\end{enumerate}
\end{enumerate}
\end{corollary}

{\bf Proof: }
    We identify the rings $R_i$ with their embeddings in $R$.
    For the first part, we use that $R^\times$ is the direct product of the unit groups $R_i^\times$ 
    and apply Lemma~\ref{lma:nrunits}.
    For the second part, we denote a minimal right ideal of size $p_i^{s_i}$ in $R_i$ by $M_i$.
    $M_i$ is also minimal when considered as a right ideal in $R$.
    By the minimality assumption, $M_i = a_iR$ for some $a_i\in R$.
    For $T\subseteq\{1,\ldots,d\}$ let $a_T = \sum_{i\in T}a_i$ and $M_T$ be the right ideal 
    $\sum_{i\in T}M_i$. Then $a_T$ generates $M_T$ as a right ideal.
    Conversely, an element $y\in M_{{\{1,\ldots,d\}}}$ generates $M_T$ as a right ideal 
    if and only if $y=\sum_{i \in T}y_i$ for nonzero $y_i \in M_i$.
    Indeed $M_T= \bigoplus_{i \in T} M_i$, and so from Lemma \ref{lemsum} it is easy to check that 
    $w(a_T) = W_T$.
    Furthermore, the right ideals of $R$ contained in $M_T$ are exactly the right ideals 
    $M_V$ with $V\subseteq T$.

    Since $C$ is regular, there is a codeword $c\in C$ with a unit entry in some coordinate.
    For $T\subseteq\{1,\ldots,d\}$ we set $c_T = a_T c$ and define non-negative integers $x_T$ as the number of elements of $c_{{\{1,\ldots,d\}}}$ that generate the right ideal $M_T$ of $R$.
    All the components of $c_{{\{1,\ldots,d\}}}$ are contained in $M_{{\{1,\ldots,d\}}}$, so condition~(i) holds true.
    It remains to show that condition~(ii) is satisfied for any $T\subseteq\{1,\ldots,d\}$, $T\neq\emptyset$:
    For $U\subseteq T$, the number of elements in $c_T$ generating the the right ideal $M_U$ is
\[
L_U = \sum_{{V\subseteq\{1,\ldots,d\}}\atop{V\cap T = U}} x_V\text{.}
\]
    Since $c$ contains a unit and $T\neq\emptyset$, the codeword $c_T$ is not the zero word.
    Using the fact that $C$ is proper, the homogeneous weight of $c_T$ is
\[
    w(c_T) = \sum_{U\subseteq T} L_U W_U\in\{w_1,w_2\}\text{.}
\]
\qed

\begin{remark}
Every finite ring is a direct product of rings of pairwise coprime prime power order.
We call this its \emph{primary decomposition}.
Corollary~\ref{coralext} can always be applied to the primary decomposition of $R$.
\end{remark}




\section{Computer Search}
In this section, a `two-weight code' will mean a primitive, proper, regular, projective two-weight code. 

Corollaries \ref{corweight}, \ref{cormulti} and \ref{coralext} can be used to analyse tables of feasible parameters of strongly regular graphs to see which might arise from two-weight codes over rings.  Feasible parameter sets for graphs of order at most $1300$ can be read at \cite{b}. 
In the following, we consider a feasible parameter set $(N,k,\lambda,\mu)$ of a strongly regular graph $G$ with restricted eigenvalues $\rho_1 < \rho_2$ of multiplicities $m_1$ and $m_2$, respectively.
If $G = \Gamma(C)$ for a proper, regular projective code $C$ of order $N$, length $n$ and nonzero weights $w_1 < w_2$ of respective frequencies $k$ and $N - k - 1$ then 
\begin{equation}\label{eqint1}
\theta := \frac{N}{\rho_2-\rho_1} \in \Z.
\end{equation} 

Moreover, for $\gamma = \abs{R^\times}$ we have $w_1=-(\rho_1+1) \theta$ and $w_2 = -\rho_1 \theta$.
Let $N=\abs{C}=\prod_{i=1}^d p_i^{t_i}$ for some distinct primes $p_i$ and positive integers $t_i$.
We apply Corollary~\ref{coralext} to the primary decomposition of $R$:
The values $\abs{R^\times}$ and $s_1,\ldots,s_d$ satisfy 
\begin{eqnarray}
&\label{eqintdiv1}1\leq s_i \leq t_i \text{ for all }i\in\{1,\ldots,d\}\text{,}\\
&\label{eqintdiv2}\prod_{i=1}^d (p_i^{s_i}-1)\text{ divides } \abs{R^\times} \text{ divides } m_2\quad\text{and}\\
&\label{eqintdiv3}\text{There exists a solution $(x_T)_{T\subseteq{\{1,\ldots,d\}}}$ in Corollary~\ref{coralext}(b).}
\end{eqnarray}
Note that these constraints allow extensive searching through the listings of feasible parameter sets without further specifying the underlying ring. 

In the remaining cases we may fix a candidate ring satisfying the above and search for a generator matrix of a two-weight code of length $n=\frac{m_2}{\abs{R^\times}}$.
The search algorithm for the required generator matrices models the problem in terms of a system of Diophantine equations.


\subsection{Choice of $R$}
In order to conduct a complete search for a given feasible parameter set $(N,k,\lambda,\mu)$ of a strongly regular graph with restricted eigenvalues $\rho_1^{m_1},\rho_2^{m_2}$ we first must have a classification of all possible rings $R$ over which there exists an $R$-module of size $N$.
Using the primary decomposition of $R$, this amounts to a classification problem on rings of prime power order.
For rings of order $p,p^2,p^3$ for a prime $p$ a complete classification can be read in \cite{r}; those with the Frobenius property are listed below.

\begin{lemma}\label{thp3}
Up to isomorphism, all finite Frobenius rings of order $p$, $p^2$ and $p^3$, $p$ prime, are direct products of chain rings.
They are given as follows ($s$ denotes the size $p^s$ of the minimal right ideals):
\[
\begin{array}{|ccccc|}
\hline
\text{order} & R                                 & \abs{ R^\times}    & s   & (R,+)\\
\hline
p            & \Z_p                              & p-1                & 1   & \Z_p \\
\hline
p^2          & \Z_p[X]/(X^2)                     & p(p-1)             & 1   & \Z_p\times\Z_p \\
p^2          & \GF(p^2)                          & (p+1)(p-1)         & 2   & \Z_p\times\Z_p \\
p^2          & \Z_p \times \Z_p                  & (p-1)^2            & 1   & \Z_p\times\Z_p \\
p^2          & \Z_{p^2}                          & p(p-1)             & 1   & \Z_{p^2} \\
\hline
p^3          & \GF(p^3)                          & (p^2 + p + 1)(p-1) & 3   & \Z_p\times\Z_p\times\Z_p \\
p^3          & \Z_p[X]/(X^3)                     & p^2(p-1)           & 1   & \Z_p\times\Z_p\times\Z_p \\
p^3          & \Z_p\times\Z_p\times\Z_p          & (p-1)^3            & 1   & \Z_p\times\Z_p\times\Z_p \\
p^3          & \Z_p[X]/(X^2)\times\Z_p           & p(p-1)^2           & 1   & \Z_p\times\Z_p\times\Z_p \\
p^3          & \GF(p^2)\times\Z_p                & (p+1)(p-1)^2       & 1,2 & \Z_p\times\Z_p\times\Z_p \\
p^3          & \Z_{p^2}[X]/(X^2 - p, X^3)        & p^2(p-1)           & 1   & \Z_{p^2}\times\Z_p \\
p^3          & \Z_{p^2}[X]/(X^2 - \alpha p, X^3) & p^2(p-1)           & 1   & \Z_{p^2}\times\Z_p \\
p^3          & \Z_{p^2} \times \Z_p              & p(p-1)^2           & 1   & \Z_{p^2}\times\Z_p \\
p^3          & \Z_{p^3}                          & p^2(p-1)           & 1   & \Z_{p^3} \\
\hline
\end{array}
\]
where for the ring $\Z_{p^2}[X]/(X^2 - \alpha p, X^3)$, $p$ is odd and $\alpha \in \Z_{p^2}$ is not a square modulo $p$.
\end{lemma}

The finite Frobenius rings of order $p^4$ and $p^5$ can be derived from \cite{dop} and \cite{cw1, cw2}, respectively; the classification for such rings is much more complex.
A general classification for arbitrary prime power order is not known.


\subsection{Choice of the isomorphism type of $C$ as an $R$-module}
\label{subsect:shape}
Once $R$ is fixed, a complete feasible search requires knowledge of the structure of submodules of $R^n$ in terms of a canonical information set and generator matrix. Therefore, we restrict the search to modules over direct products of chain rings; specifically we restrict the search to modules over rings of order not divisible by a fourth prime power.\footnote{The order of a finite Frobenius ring that is not a direct product of finite chain rings necessarily is divisible by the fourth power of some prime, for example, the matrix ring $\Z_p^{2\times 2}$.} This still allows us to address existence of two-weight codes of most orders less than $1300$. In addition, we exclude the case when $R$ is isomorphic to a finite field from the search as this is a return to the classical case.
A finite chain ring is one whose ideals are linearly ordered. Such rings have prime power order. See \cite{cd,nec,Nechaev-FiniteRings-2008} for more on the theory of finite chain rings. For the remainder of this section we assume that $R$ is a direct product of chain rings $R_1,\ldots,R_s$ where for each $i$, $R_i$ has residue field $\GF(q_i)$, chain length $\ell_i$ and ideal chain ${0} = R_i\theta_i^{\ell_i} \subset R_i\theta_i^{\ell_i-1} \subset \ldots\subset R_i\theta_i = \rad R_i \subset R_i$, for some $\theta_i \in R_i$. 
For each pair $i,j$, let $I(i,j) =R_i \theta^{\ell_i - j}$ and let $S(i,j)$ be a transversal of $R_i / I(i,j)$ in $R_i$.

Let $\pi_i : R\rightarrow R$ be the projection to the $i$-th component of $R$, that is
\[
(\pi_i(a_1,\ldots,a_s))_j = \begin{cases}a_i & \text{if }i=j\text{,}\\0 & \text{otherwise.}\end{cases}
\]
$\pi_i$ canonically extends to mappings on vectors over $R$. In particular every left $R$-module $M$ has a unique expression in the form
$M = \bigoplus_{i=1}^s \pi_i(M)$.


The following follows directly from \cite[Theorem 2.2]{Honold-Landjev-2000-EJC7:R11}. It describes the form of a generator matrix of a module over a direct sum of chain rings. As usual, we write $\lambda \vdash t$ if $\lambda$ is a partition of $t$. By $\ell(\mu)$ we denote the \emph{length} of a partition $\mu$, that is, the number of its non-zero summands.

\begin{lemma}
\label{lma:modulestructure}
Let $M < {_R}R^n$. Then $M$ is a direct sum of cyclic left $R$-modules.
There is a unique $s$-tuple $\Lambda = (\lambda_1,\ldots,\lambda_s)$ of partitions $\Lambda_i = (\lambda^i_1,\lambda^i_2,\ldots \lambda^i_{e_i})\vdash \log _{q_i} \abs{\pi_i(M)}$ such that
\[
M \cong \bigoplus_{i = 1}^s \bigoplus_{j = 1}^{e_i} I(i,\lambda^i_j).
\]
$M$ is the row space of an $(e_1 + \ldots + e_s) \times n$ matrix $Y$ over $R$, whose columns are elements of $\bigoplus_{i=1}^s\bigoplus_{j=1}^{e_i}I(i,\lambda^i_j)$.
Each element of $M$ can be uniquely written in the form $xY$, where $x$ is a row vector in
\[
\bigtimes_{i = 1}^s \bigtimes_{j = 1}^{e_i} S(i,\lambda_j^i)\text{.}
\]
\end{lemma}

The $\Lambda$ of Lemma \ref{lma:modulestructure} is called the {\em shape} of the module $M$ over $R$. 
To generate all possible isomorphism types of $C$ as an $R$-module, we loop over all shapes $\Lambda$ given by Lemma~\ref{lma:modulestructure} leading to an $R$-module of size $N$.
Knowing the shape of a code means we may assume that $C$ is given by a generator matrix $Y$ as in Lemma \ref{lma:modulestructure}.
Thus only certain vectors of $R^n$ may appear as columns in $Y$. 
There are further restrictions we can impose on $\Lambda$.
Since $C$ is projective, we have $\lambda_1^i = \ell_i$ for all $i$.
Furthermore if $R_i = R_j$ for some $i\neq j$, the ring $R$ has an automorphism interchanging the $i$-th and the $j$-th component, so up to code isomorphism we may assume $\lambda_i \leq \lambda_j$ with respect to some fixed linear ordering on the partitions.





\subsection{Diophantine equation system}
\label{subsect:dio_eq}

Once the base ring $R$ and the shape $\Lambda = (\lambda_1,\ldots,\lambda_s)$ with partitions $\lambda_i = (\lambda_1^i,\ldots,\lambda_{e_i}^i)$ is fixed, the existence problem for the two-weight code $C$ can be reformulated as a Diophantine equation system.
This is done by adapting the method of \cite{Kohnert-2007-DAM155[1]:1451-1457}.

By Lemma~\ref{lma:modulestructure}, $C$ has a generator matrix $Y \in R^{(e_1 + \ldots + e_s)\times n}$ whose columns are elements of the module $M = \bigoplus_{i = 1}^s \bigoplus_{j = 1}^{e_i} I(i,\lambda_i^j)$, and each codeword of $C$ can be uniquely written in the form $xY$, where $x$ is a row vector in $X = \bigtimes_{i = 1}^s \bigtimes_{j = 1}^{e_i} S(i,\lambda_i^j)$.
Let ${\cal O}$ be a set of right-projective representatives of the regular vectors in $M$. Let $r=\abs{\cal O}$ and let $t=\abs{X}$.
Let $v \in \{0,1\}^r$ be the characteristic vector of the columns of $Y$ in ${\cal O}$ whose entries are labelled by the elements of ${\cal O}$.
That is, for each $y \in {\cal O}$, $v_y = 1$ if $y$ is right projectively equivalent to a column of $Y$ and is zero otherwise. Since $C$ is regular and projective, each column of $Y$ corresponds to exactly one element of ${\cal O}$. In particular, any two codes determine the same characteristic vector if and only if they are monomially equivalent.

Define a matrix $W \in \Z^{(t - 1)\times r}$ whose rows are labelled by the elements of $X \setminus \{0\}$, whose columns are labelled by the elements of ${\cal O}$ and which satisfies $W_{x,y} = w(\langle x,y \rangle) = w((xM)_y)$ for each $x \in X$ and $y \in {\cal O}$.
The weight of the codeword $xY$ is given by $w(xY) = (Wv)_x$.
Let $0<w_1<w_2 \in \Z$.
Let $Z$ be the $(t-1) \times (t-1)$ diagonal matrix $(w_2 - w_1)I$.
Then $C$ is a two-weight code of the given parameters if and only if there is a vector $u\in \{0,1\}^{t - 1}$ of auxiliary variables such that $(v,u)\in\{0,1\}^{r + t - 1}$ is a solution of the Diophantine equation system
\begin{equation}\label{eqdio}
\left(\begin{array}{cc}
W& Z \\
1\ldots 1 & 0\ldots 0\\
0\ldots 0 & 1\ldots 1\\
\end{array}\right)
\left(\begin{array}{c} v \\  u\end{array}\right) =
\left(\begin{array}{c} w_2 \\ \vdots \\ w_2 \\ n \\ k
\end{array}\right).
\end{equation}

Therefore, the existence of $C$ can be decided by applying an integer linear program solver to this equation system.
In the case that a solution exists, the solution part $v$ gives us the columns of a generator matrix of a suitable code $C$.
The solver we used is described in \cite{wa1,wa2}.

\subsection{Results}

We performed a search for two-weight codes $C$ of order at most $1300$ as follows:
\begin{enumerate}
\item
   We listed all parameter sets $(N,k,\lambda,\mu),\rho_1^{m_1},\rho_2^{m_2}$ given in \cite{b} satisfying the integrality condition (\ref{eqint1}).
\item
   For each such paramter set in the above list, we write $N = \prod_{i=1}^d p_i^{t_i}$ with $p_i$ pairwise distinct primes and assume that the underlying ring has the primary decomposition $R = R_1\oplus\ldots\oplus R_d$ with $\abs{R_i} = p_i^{r_i}$ and $1\leq r_i\leq t_i$ (since $C$ is projective and regular over $R$).
We listed all possibilities for the tuples $(\abs{R^\times},s_1,\ldots,s_d)$ satisfying the conditions (\ref{eqintdiv1}), (\ref{eqintdiv2}) and (\ref{eqintdiv3}).
\item
   For each tuple $(\abs{R^\times},s_1,\ldots,s_d)$ in the above list, using Lemma~\ref{thp3} we generated all finite Frobenius rings $R = R_1\oplus\ldots\oplus R_d$ with $\abs{R_i} = p_i^{r_i}$ such that 
   \begin{itemize}
      \item For all $i\in\{1,\ldots,d\}$, $1\leq r_i \leq s_i$.
      \item
         $R$ has the prescribed number of units.
      \item
         For all $i\in\{1,\ldots,d\}$, $R_i$ has a minimal right ideal of size $p_i^{s_i}$.
      \item
         $R$ is not a finite field and $\abs{R}$ is not divisible by a fourth prime power. 
   \end{itemize}  
\item
   For each such ring $R$ we listed all the potential shapes $\Lambda$ of $C$ as discussed in Section~\ref{subsect:shape}.
\item
   For each such shape $\Lambda$ we sought a solution to a system of Diophantine equations as in (\ref{eqdio}).    
\end{enumerate}


In every case that a two-weight code was found, $N$ was the square of a prime power.\footnote{There are classical two-weight codes of dimension $9$ over $\Z_2$, which suggests that the square property may not be true for larger orders.} 
The computations were completed successfully 
for all orders $N$ except
\[
256, 324, 486, 512, 576, 640, 729, 768, 800, 1024, 1296\text{.}
\]
For each of the above orders, there is some parameter set that yields an equation system too large for the solver to terminate within a reasonable time limit.

If $N = p_1^4 p_2 \ldots p_t$ for distinct primes $p_i$ and all computations were successfully completed under the assumption that $|R|$ is not divisible by a fourth power, then we have the 
following additional existence criteria: If there exists a two-weight $C$ over $R$ such that $\Gamma(C)$ is an $(N,K,\lambda,\mu)$ strongly regular graph then $N = \abs{C} = \abs{R}$ and hence $n = 1$ as $C$ is projective.
Then the equation system in Corollary~\ref{corweight}(b) must have a solution for $n = 1$.
This eliminates the putative parameter sets $( 162, 92, 46, 60 ),( 162, 138, 117, 120 )$ and $( 162, 140, 121, 120 )$ as arising from a two-weight code.

The only orders $N \leq 1300$, not the square of a prime power, for which the existence of a two-weight code over some finite Frobenius ring has not been ruled out by our results are:
\[
96, 144, 243, 288, 320, 324, 400, 486, 512, 576, 640, 768, 784, 800, 1200, 1296\text{.}
\]


There are currently $2960$ (up to complements: $1514$) feasible parameter sets for strongly regular graphs on at most $1300$ vertices for which the actual existence of a corresponding graph is not yet known \cite{b}.
After applying (\ref{eqint1}), there remain $867$ such parameter sets.
Filtering further for conditions (\ref{eqintdiv1}), (\ref{eqintdiv2}) and (\ref{eqintdiv3}), we are left with $129$ parameter sets.
In the end, there are only $82$ cases where it is not known whether or not such a graph may arise from a two-weight code. We remark that this existence question is not invariant under taking complements. Among the $82$ cases there are no self-complementary parameter sets (recall we must have $\rho_1<\rho_2$), and exactly $6$ pairs of parameter sets complementary to each other, these are:
\[
\begin{array}{ccc}
(729,140,13,30) & \text{and} & (729,588,477,462)\text{,} \\
(1024,165,8,30) & \text{and} & (1024,858,722,702)\text{,} \\
(1024,363,122,132) & \text{and} & (1024,660,428,420)\text{,} \\
(1024,396,148,156) & \text{and} & (1024,627,386,380)\text{,} \\
(1024,429,176,182) & \text{and} & (1024,594,346,342)\text{,} \\
(1024,462,206,210) & \text{and} & (1024,561,308,306)\text{.}
\end{array}
\]
For the remaining $70$ parameter sets the complement has been eliminated.
We list these $82$ cases in detail below:

\begin{longtable}{|cccc|cc|cc|}
\hline
$N$ & $k$ & $\lambda$ & $\mu$ & $\rho_1^{m_1}$ & $\rho_2^{m_2}$ & \multicolumn{2}{|c|}{weights ($\gamma = \abs{R^\times}$)} \\
\hline
\endhead
\hline
\endfoot
$96$ & $45$ & $24$ & $18$ & $-3^{75}$ & $9^{20}$ & $16^{45}$ & $24^{50}$ \\
$144$ & $91$ & $58$ & $56$ & $-5^{91}$ & $7^{52}$ & $48^{91}$ & $60^{52}$ \\
$288$ & $41$ & $4$ & $6$ & $-7^{123}$ & $5^{164}$ & $144^{41}$ & $168^{246}$ \\
$288$ & $42$ & $6$ & $6$ & $-6^{147}$ & $6^{140}$ & $120^{42}$ & $144^{245}$ \\
$288$ & $123$ & $42$ & $60$ & $-21^{41}$ & $3^{246}$ & $240^{123}$ & $252^{164}$ \\
$288$ & $140$ & $76$ & $60$ & $-4^{245}$ & $20^{42}$ & $36^{140}$ & $48^{147}$ \\
$320$ & $145$ & $60$ & $70$ & $-15^{87}$ & $5^{232}$ & $224^{145}$ & $240^{174}$ \\
$320$ & $154$ & $78$ & $70$ & $-6^{231}$ & $14^{88}$ & $80^{154}$ & $96^{165}$ \\
$324$ & $95$ & $22$ & $30$ & $-13^{95}$ & $5^{228}$ & $216^{95}$ & $234^{228}$ \\
$324$ & $102$ & $36$ & $30$ & $-6^{221}$ & $12^{102}$ & $90^{102}$ & $108^{221}$ \\
$324$ & $133$ & $52$ & $56$ & $-11^{133}$ & $7^{190}$ & $180^{133}$ & $198^{190}$ \\
$324$ & $136$ & $58$ & $56$ & $-8^{187}$ & $10^{136}$ & $126^{136}$ & $144^{187}$ \\
$324$ & $204$ & $126$ & $132$ & $-12^{119}$ & $6^{204}$ & $198^{204}$ & $216^{119}$ \\
$324$ & $209$ & $136$ & $132$ & $-7^{209}$ & $11^{114}$ & $108^{209}$ & $126^{114}$ \\
$324$ & $247$ & $190$ & $182$ & $-5^{247}$ & $13^{76}$ & $72^{247}$ & $90^{76}$ \\
$400$ & $147$ & $50$ & $56$ & $-13^{147}$ & $7^{252}$ & $240^{147}$ & $260^{252}$ \\
$400$ & $152$ & $60$ & $56$ & $-8^{247}$ & $12^{152}$ & $140^{152}$ & $160^{247}$ \\
$400$ & $228$ & $128$ & $132$ & $-12^{171}$ & $8^{228}$ & $220^{228}$ & $240^{171}$ \\
$400$ & $231$ & $134$ & $132$ & $-9^{231}$ & $11^{168}$ & $160^{231}$ & $180^{168}$ \\
$400$ & $315$ & $250$ & $240$ & $-5^{315}$ & $15^{84}$ & $80^{315}$ & $100^{84}$ \\
$486$ & $100$ & $22$ & $20$ & $-8^{275}$ & $10^{210}$ & $189^{100}$ & $216^{385}$ \\
$486$ & $194$ & $67$ & $84$ & $-22^{97}$ & $5^{388}$ & $378^{194}$ & $396^{291}$ \\
$486$ & $210$ & $99$ & $84$ & $-6^{385}$ & $21^{100}$ & $90^{210}$ & $108^{275}$ \\
$486$ & $388$ & $310$ & $308$ & $-8^{291}$ & $10^{194}$ & $189^{388}$ & $216^{97}$ \\
$576$ & $125$ & $16$ & $30$ & $-19^{125}$ & $5^{450}$ & $432^{125}$ & $456^{450}$ \\
$576$ & $175$ & $46$ & $56$ & $-17^{175}$ & $7^{400}$ & $384^{175}$ & $408^{400}$ \\
$576$ & $225$ & $84$ & $90$ & $-15^{225}$ & $9^{350}$ & $336^{225}$ & $360^{350}$ \\
$576$ & $230$ & $94$ & $90$ & $-10^{345}$ & $14^{230}$ & $216^{230}$ & $240^{345}$ \\
$576$ & $322$ & $178$ & $182$ & $-14^{253}$ & $10^{322}$ & $312^{322}$ & $336^{253}$ \\
$576$ & $325$ & $184$ & $182$ & $-11^{325}$ & $13^{250}$ & $240^{325}$ & $264^{250}$ \\
$576$ & $375$ & $246$ & $240$ & $-9^{375}$ & $15^{200}$ & $192^{375}$ & $216^{200}$ \\
$576$ & $425$ & $316$ & $306$ & $-7^{425}$ & $17^{150}$ & $144^{425}$ & $168^{150}$ \\
$576$ & $475$ & $394$ & $380$ & $-5^{475}$ & $19^{100}$ & $96^{475}$ & $120^{100}$ \\
$640$ & $284$ & $108$ & $140$ & $-36^{71}$ & $4^{568}$ & $560^{284}$ & $576^{355}$ \\
$640$ & $315$ & $170$ & $140$ & $-5^{567}$ & $35^{72}$ & $64^{315}$ & $80^{324}$ \\
$729$ & $140$ & $13$ & $30$ & $-22^{140}$ & $5^{588}$ & $567^{140}$ & $594^{588}$ \\
$729$ & $588$ & $477$ & $462$ & $-6^{588}$ & $21^{140}$ & $135^{588}$ & $162^{140}$ \\
$768$ & $59$ & $10$ & $4$ & $-5^{531}$ & $11^{236}$ & $192^{59}$ & $240^{708}$ \\
$768$ & $177$ & $36$ & $42$ & $-15^{295}$ & $9^{472}$ & $448^{177}$ & $480^{590}$ \\
$768$ & $182$ & $46$ & $42$ & $-10^{455}$ & $14^{312}$ & $288^{182}$ & $320^{585}$ \\
$768$ & $295$ & $102$ & $120$ & $-25^{177}$ & $7^{590}$ & $576^{295}$ & $600^{472}$ \\
$768$ & $312$ & $136$ & $120$ & $-8^{585}$ & $24^{182}$ & $168^{312}$ & $192^{455}$ \\
$768$ & $531$ & $354$ & $396$ & $-45^{59}$ & $3^{708}$ & $704^{531}$ & $720^{236}$ \\
$784$ & $261$ & $80$ & $90$ & $-19^{261}$ & $9^{522}$ & $504^{261}$ & $532^{522}$ \\
$784$ & $324$ & $136$ & $132$ & $-12^{459}$ & $16^{324}$ & $308^{324}$ & $336^{459}$ \\
$784$ & $432$ & $236$ & $240$ & $-16^{351}$ & $12^{432}$ & $420^{432}$ & $448^{351}$ \\
$784$ & $435$ & $242$ & $240$ & $-13^{435}$ & $15^{348}$ & $336^{435}$ & $364^{348}$ \\
$784$ & $609$ & $476$ & $462$ & $-7^{609}$ & $21^{174}$ & $168^{609}$ & $196^{174}$ \\
$800$ & $705$ & $620$ & $630$ & $-15^{235}$ & $5^{564}$ & $560^{705}$ & $600^{94}$ \\
$800$ & $714$ & $638$ & $630$ & $-6^{595}$ & $14^{204}$ & $200^{714}$ & $240^{85}$ \\
$1024$ & $165$ & $8$ & $30$ & $-27^{165}$ & $5^{858}$ & $832^{165}$ & $864^{858}$ \\
$1024$ & $363$ & $122$ & $132$ & $-21^{363}$ & $11^{660}$ & $640^{363}$ & $672^{660}$ \\
$1024$ & $396$ & $148$ & $156$ & $-20^{396}$ & $12^{627}$ & $608^{396}$ & $640^{627}$ \\
$1024$ & $429$ & $176$ & $182$ & $-19^{429}$ & $13^{594}$ & $576^{429}$ & $608^{594}$ \\
$1024$ & $462$ & $206$ & $210$ & $-18^{462}$ & $14^{561}$ & $544^{462}$ & $576^{561}$ \\
$1024$ & $561$ & $308$ & $306$ & $-15^{561}$ & $17^{462}$ & $448^{561}$ & $480^{462}$ \\
$1024$ & $594$ & $346$ & $342$ & $-14^{594}$ & $18^{429}$ & $416^{594}$ & $448^{429}$ \\
$1024$ & $627$ & $386$ & $380$ & $-13^{627}$ & $19^{396}$ & $384^{627}$ & $416^{396}$ \\
$1024$ & $660$ & $428$ & $420$ & $-12^{660}$ & $20^{363}$ & $352^{660}$ & $384^{363}$ \\
$1024$ & $858$ & $722$ & $702$ & $-6^{858}$ & $26^{165}$ & $160^{858}$ & $192^{165}$ \\
$1200$ & $218$ & $28$ & $42$ & $-22^{327}$ & $8^{872}$ & $840^{218}$ & $880^{981}$ \\
$1200$ & $231$ & $54$ & $42$ & $-9^{847}$ & $21^{352}$ & $320^{231}$ & $360^{968}$ \\
$1200$ & $847$ & $590$ & $616$ & $-33^{231}$ & $7^{968}$ & $960^{847}$ & $990^{352}$ \\
$1296$ & $185$ & $4$ & $30$ & $-31^{185}$ & $5^{1110}$ & $1080^{185}$ & $1116^{1110}$ \\
$1296$ & $259$ & $34$ & $56$ & $-29^{259}$ & $7^{1036}$ & $1008^{259}$ & $1044^{1036}$ \\
$1296$ & $333$ & $72$ & $90$ & $-27^{333}$ & $9^{962}$ & $936^{333}$ & $972^{962}$ \\
$1296$ & $407$ & $118$ & $132$ & $-25^{407}$ & $11^{888}$ & $864^{407}$ & $900^{888}$ \\
$1296$ & $420$ & $144$ & $132$ & $-12^{875}$ & $24^{420}$ & $396^{420}$ & $432^{875}$ \\
$1296$ & $481$ & $172$ & $182$ & $-23^{481}$ & $13^{814}$ & $792^{481}$ & $828^{814}$ \\
$1296$ & $490$ & $190$ & $182$ & $-14^{805}$ & $22^{490}$ & $468^{490}$ & $504^{805}$ \\
$1296$ & $555$ & $234$ & $240$ & $-21^{555}$ & $15^{740}$ & $720^{555}$ & $756^{740}$ \\
$1296$ & $560$ & $244$ & $240$ & $-16^{735}$ & $20^{560}$ & $540^{560}$ & $576^{735}$ \\
$1296$ & $700$ & $376$ & $380$ & $-20^{595}$ & $16^{700}$ & $684^{700}$ & $720^{595}$ \\
$1296$ & $703$ & $382$ & $380$ & $-17^{703}$ & $19^{592}$ & $576^{703}$ & $612^{592}$ \\
$1296$ & $770$ & $454$ & $462$ & $-22^{525}$ & $14^{770}$ & $756^{770}$ & $792^{525}$ \\
$1296$ & $777$ & $468$ & $462$ & $-15^{777}$ & $21^{518}$ & $504^{777}$ & $540^{518}$ \\
$1296$ & $840$ & $540$ & $552$ & $-24^{455}$ & $12^{840}$ & $828^{840}$ & $864^{455}$ \\
$1296$ & $851$ & $562$ & $552$ & $-13^{851}$ & $23^{444}$ & $432^{851}$ & $468^{444}$ \\
$1296$ & $910$ & $634$ & $650$ & $-26^{385}$ & $10^{910}$ & $900^{910}$ & $936^{385}$ \\
$1296$ & $925$ & $664$ & $650$ & $-11^{925}$ & $25^{370}$ & $360^{925}$ & $396^{370}$ \\
$1296$ & $999$ & $774$ & $756$ & $-9^{999}$ & $27^{296}$ & $288^{999}$ & $324^{296}$ \\
$1296$ & $1073$ & $892$ & $870$ & $-7^{1073}$ & $29^{222}$ & $216^{1073}$ & $252^{222}$
\end{longtable}


Besides these nonexistence results, in some cases we also got new existence results.
For example, let $R = \GF(4)\times\GF(2)$ and consider the $R$-linear code $C$ of length $3$ given by the row space of the matrix
\[
\left[\begin{array}{ccc}
(1,1) & (1,1) & (1,0) \\
(0,0) & (0,1) & (0,1)
\end{array}\right]
.
\]
$C$ is a primitive proper regular projective two-weight code with homogeneous weight enumerator $0^1 8^9 12^6$ (scaling factor $\gamma$ set to $\abs{R^\times} = 3$).
As an $R$-module, $C$ is isomorphic to $R(1,0) \oplus R(0,1) \oplus R(0,1)$, so $C$ is not free.
The strongly regular graph generated by $C$ is isomorphic to the $4^2$-graph.
In fact, it turns out that the code $C$ is a member of an infinite family of two-weight codes codes over base rings of the form $\GF(p^k) \times \GF(p^l)$, which generate strongly
regular graphs having the parameters of those determined by orthogonal arrays of strength 2.
Details will be published in a forthcoming paper.

\section{Gray Isometries}

An immediate consequence of Corollary \ref{corgamma} is given by the following.

\begin{corollary}\label{corprimepower}
   Let $R$ have prime power order $p^s$ for some prime $p$ and positive integer $s$. Let $C$ be 
   a proper, regular,   
   projective two-weight code over $R$ with nonzero weights $w_1<w_2$, where the weight function
   is computed for $\gamma = \abs{R^\times}$. Then there exist integers
   $r$ and $t$, $r \geq 0, t > 0$, satisfying $w_1=p^r t$ and $w_2 = p^r(t+1)$.
\end{corollary}

One question that arises from Corollary \ref{corprimepower} concerns whether or not
a two-weight code of prime power order yields a graph isomorphic to one arising from a finite field. 

Let $C_1,C_2$ be a pair of two-weight codes over finite rings $R_1,R_2$ respectively, with respect to a pair of (possibly distinct) weight functions $w^1,w^2$. Let $\Gamma_i:=\Gamma(C_i)$ for $i =1,2$. Then clearly $\Gamma_1$ and $\Gamma_2$ are isomorphic graphs if there is an isometry (or scaled isometry) $\iota : (C_1,w^1) \longrightarrow (C_2,w^2)$.

We now consider the possibility that some of the known constructions of linear two-weight codes over a finite field are images of linear codes over a finite chain ring under the Gray isometry.
A number of authors have looked into extending the standard Gray isometry between $(\Z_4,w_{\rm {Lee}})$ and $(\Z_2^2,w_{{\rm Ham}})$ for the case of a finite chain ring (c.f. \cite{gs99,hn,hl,nk}). The Gray map for more general rings has been considered in \cite{hn}. If $R$ is a finite chain ring of length $n$ and residue field $\GF(q)$ there is an isometric embedding of $R$ for the homogeneous weight into $\GF(q)^{q^{n-1}}$ for the Hamming weight (which of course is homogeneous over $\GF(q)$), in which case the image of $R$ is the generalized Reed-Muller code GRM$(1,n-1)$. 

For example, in \cite[Prop 6.2]{bgh} a construction is given for a two-weight code $C$ over a finite chain ring $R$ of length 2 and having residue field $\GF(q)$. $C$ has a $2 \times s(q+1)$ generator matrix whose $s(q+1)$ columns comprise $s$ distinct elements from each equivalence class of $q$ points in the projective Hjelmslev line over $R$, $1 \leq s \leq q-1$. 
Then $C$ has $q^4$ codewords and non-zero homogeneous weights
$$w_1 = q^2(qs-1) \mbox{ and } w_2 = q^3s = q^2(qs),$$
for $\gamma = \abs{R^\times} = q^2-q$. We easily solve for $k,\rho_1,\rho_2$ to find $k =s(q^3-q) $, $\rho_1 = -qs$ and $\rho_2 = q^2 - qs$, from which we may conclude, using Lemma \ref{lemsrg}, that $\Gamma(C)$ is a strongly regular graph with 
parameters 
$$(q^4,s(q^3-q),\lambda = q^2(s^2+1)-3qs,\mu = qs(qs-1)).$$
For the case $s=1$, $R=\Z_4,\Z_9$ and $\GR(4,2)$ (the Galois ring of order $16$ and characteristic $4$) the corresponding two-weight codes have lengths $3,4$ and $5$, sizes $16,81$ and $256$ and yield strongly regular graphs with parameters $(16,4,2,2),(81,24,9,6)$ and $(256,60,20,12)$, respectively. While each such graph has the same parameters as the Cayley graph of a binary $[6,4,2,4]$ two-weight code, a ternary $[12,4,6,9]$ two-weight code, and a $\GF(4)$-$[20,4,12,16]$ two-weight code, respectively, we have verified by computer search that no Gray image of any of the $4,77$ or $1023$ distinct $\Z_4,\Z_9$ or $\GR(4,2)$  codes constructed as above is $\Z_2,\Z_3$ or $\GF(4)$-linear. Obviously this is only a partial result towards the question of graph isomorphism for the codes considered here.

\end{document}